\documentclass{amsproc}
  
\usepackage{graphicx}

\newtheorem{thm}{Theorem}[section]
\newtheorem{prop}[thm]{Proposition}
\newtheorem{lemma}[thm]{Lemma}

\theoremstyle{definition}
\newtheorem{defn}[thm]{Definition}
\newtheorem{ex}[thm]{Example}

\theoremstyle{remark}
\newtheorem{remark}[thm]{Remark}

\copyrightinfo{2001}
    {American Mathematical Society}

\numberwithin{equation}{section}

\def\R{\mathbb{R}}

\def\P{\mathbb{P}}

\begin{document}

\title{Special Isothermic Surfaces and Solitons}

\author{Emilio Musso}
\address{Dipartimento di Matematica Pura ed Applicata,
Universit\`a degli Studi di L'Aquila, Via Vetoio, I-67100 L'Aquila, Italy}
\email{musso@univaq.it}

\author{Lorenzo Nicolodi}
\address{Di\-par\-ti\-men\-to di Ma\-te\-ma\-ti\-ca, 
Uni\-ver\-si\-t\`a degli Studi di Parma, Via M. D'Azeglio 85/A, 
I-43100 Parma, Italy}
\email{nicolodi@mat.uniroma1.it}

\thanks{This research was partially supported by the MURST project 
\textit{Propriet\`a Geometriche delle Variet\`a Reali e Complesse}, 
by the GNSAGA of the CNR, and by the European Contract Human Potential 
Programme, Research Training Network HPRN-CT-2000-00101 (EDGE).
Part of this work was done while the first author was visiting the 
Laboratoire de Math\'ematiques J.-A. Dieudonn\'e, 
Universit\'e de Nice Sophia-Antipolis. We thank these institutions for 
their support and/or hospitality.}

\date{April 9, 2001}



\begin{abstract} 
We establish a correspondence between
Darboux's special isothermic surfaces of type $(A,0,C,D)$ and the
solutions of the second order p.d.e. 
$\Phi\Delta\Phi -|\nabla\Phi|^{2} +\Phi^{4}
=s,  s\in \R$. We then use the classical Darboux transformation
for isothermic surfaces to construct a B\"acklund transformation
for this equation and prove a superposition formula for its
solutions. As an application we discuss 1 and 2-soliton solutions
and the corresponding surfaces.
\end{abstract}

\maketitle

\section{Introduction}\label{intro}
The theory of isothermic surfaces in conformal geometry has been
the focus for intense research over the past years due to its
relation with the theory of integrable systems; see for instance
\cite{BDPT, Bu, BHPP, CGS, Ci, HJ, HP, HMN} and the literature
therein.

\medskip
\noindent Among isothermic surfaces Darboux \cite{Da}
distinguished
the class of {\it special isothermic surfaces}. These were
introduced in connection with the problem of isometric deformation
of quadrics and have been investigated by Bianchi \cite{B1, B2}
and Calapso \cite{C2}. 
An isothermic immersion $f : U\subset
\R^2\to \R^{3}$
admits conformal curvature line coordinates $x,y$ for which its
first and second fundamental forms read 
$$
I=\Theta^{2}(dx^{2}+dy^{2}),\quad
II=\Theta^{2}(h_{1}dx^{2}+h_{2}dy^{2}), 
$$ 
where $\Theta$ is a
nowhere vanishing smooth function and $h_{1}$, $h_{2}$ denote the
principal curvatures.
The immersion $f$ is called special of type $(A,B,C,D)$ if its
mean curvature $H$ satisfies the equation
\begin{equation}
4\Theta^{2}|\nabla H|^{2}+M^{2}+2AM+2BH +2CL+D=0,
\end{equation}
where $A,B,C,D$ are real constants and $M=-HL$ being
$L=\Theta^{2}(h_{1}-h_{2})$. Examples include surfaces of constant
mean curvature (cmc).

\medskip
\noindent  This article has its origins in the two seminal papers
of Bianchi \cite{B1, B2} on isothermic surfaces and deals with special 
isothermic surfaces of type $(A,0,C,D)$, henceforth simply referred to as 
special isothermic. 
These surfaces are invariant under the group of 
conformal transformations. Since the work of Bianchi, it has been known that 
umbilic free special isothermic surfaces are conformally 
equivalent to cmc surfaces in 3-dimensional space 
forms\footnote{Observe that cmc immersions in 3-dimensional space forms are
obtained as T-transforms (spectral deformations) of minimal immersions in 
space forms (Willmore isothermic surfaces). For a recent account of these
facts see \cite{HMN}, \cite{CM}, and the monograph of Burstall
\cite{Bu}.}. Examples of special isothermic surfaces 
with umbilic lines can be constructed by revolving
elastic curves in the hyperbolic half-plane about the
boundary at infinity (cf. Section \ref{se.geom}).

\medskip
\noindent We prove that special isothermic immersions 
are in correspondence with the solutions $\Phi(x,y)$ of the
partial differential equation \footnote{For a special isothermic
immersion, the corresponding solution is the Calapso potential
$\Phi$ defined by $\Phi^2(dx^2 + dy^2)={(1/4)}(h_1-h_2)^2I$ 
which only depends on the conformal class of the induced metric.}
\begin{equation}
\Phi\Delta\Phi -|\nabla\Phi|^{2} +\Phi^{4} =s, \quad s\in
\R.\label{solitoneq}
\end{equation}
Further, we show how the geometric properties of the Darboux
transformation of special isothermic surfaces apply to obtain
analytic results for this equation. As for the study of Darboux
transforms of special isothermic surfaces, our work is related to
that of Hertrich-Jeromin--Pedit \cite{HP} which discusses the case
of cmc immersions in Euclidean space and that of
Hertrich-Jeromin--Musso--Nicolodi \cite{HMN} concerning immersions
with cmc $H$, $|H|=1$, in hyperbolic 3-space.

\medskip
\noindent The special isothermic surface equation
(\ref{solitoneq}) can be reformulated as a zero-curvature equation
\begin{equation}
 P(m)_y - Q(m)_x - [P(m),Q(m)] = 0
\end{equation}
involving an auxiliary parameter $m$. This equation expresses the
compatibility condition for a linear differential system
\begin{equation}
V_x=-P(m)V,\quad V_y=-Q(m)V \label{linearsys}
\end{equation}
in Minkowski 5-space, where the Lax pair matrices $P(m)$, $Q(m)$
take values in ${\bf{so}}(4,1)$ (cf. Section \ref{se.lax}). In
this representation, for every solution $\Phi$ of
(\ref{solitoneq}), the spectral parameter $m$ describes a family
of special isothermic surfaces which are second order conformal
deformations of each other\footnote{They amount to a special
isothermic immersion together with its T-transforms \cite{M1} 
(see also Theorem 4.6 below).}
and are obtained by integrating the loop of
${\bf{so}}(4,1)$-valued  1-forms $P(m)dx+Q(m)dy$.

\medskip
\noindent The linear system (\ref{linearsys}) have the conserved
quantities
\begin{equation}
 C_1(V) = -2v^{0}v^{4}+\sum_{a=1}^{3}(v^{a})^{2},\quad
 C_2(V) = \Phi^{2}v^{0}-mv^{3}+\frac{s}{2}\Phi^{-2}v^{4}.
\end{equation}
The classical Darboux transforms of a special isothermic immersion
with potential $\Phi$ are constructed from the solutions $V$ of
the linear system (\ref{linearsys}) satisfying $C_1(V)=0$ (cf.
\cite{B1,B2} and the recent papers \cite{Bu,Ci,BHPP,HMN, HP}).
Moreover, if $V$ satisfies the additional constraint $C_2(V)=0$,
the Darboux transforms of special isothermic surfaces are still
special isothermic. Analytically, this is equivalent to the
statement that if $\Phi$ is a solution of (\ref{solitoneq}) and if
$V$ is a solution of the linear system satisfying
$C_{1}(V)=C_{2}(V)=0$, then $(v^{3}/v^{4})\Phi$ is again a
solution which can be regarded as a B\"acklund transformation of
(\ref{solitoneq}). This description furnishes a procedure for
generating new solutions of (\ref{solitoneq}) by solving
(\ref{linearsys}).

\medskip
\noindent In Section \ref{se.superpo} we prove a superposition
formula for the solutions of (\ref{solitoneq}). Namely, if $\Phi$
is a known solution and $\Phi_1$, $\Phi_2$ are solutions
respectively generated from $\Phi$ by solutions $V$ and $W$ of
(\ref{linearsys}) corresponding to different values of the
spectral parameter, then $$ \Phi_3= \Phi - \frac{ (h-k){v}^3{w}^3}
{v^{1}w^{1}+v^{2}w^{2}+v^{3}w^{3}-v^{0}w^{4}-v^{4}w^{0}}
\left[\frac{\Phi_2 - \Phi_1}{\Phi_1\Phi_2}\right] $$ represents a
new solution of (\ref{solitoneq}). This formula, which
geometrically amounts to the Bianchi permutability theorem for
special isothermic immersions, shows that, after the first step,
the procedure for obtaining new solutions can be carried out
without the  quadratures associated with the B\"acklund
transformation; it also provides a  method for obtaining the
multisoliton solutions of equation (\ref{solitoneq}) by algebraic
means only. As an application, we compute the 1 and 2-soliton
solutions arising from the trivial solution $\Phi=1$ of
(\ref{solitoneq}) with $s=1$.

\smallskip
\noindent The B\"acklund transformation for the fourth-order equation
defining isothermic immersions --- the so-called  Calapso equation
\footnote{For more information on the Calapso equation we refer
the reader to \cite{C1, BHPP, CGS, Bu}} --- as well as its
1-solitons have been previously considered in \cite{MNcnr} and have been used
by Bernstein \cite{Be} to obtain explicit examples of non-special
isothermic tori with spherical curvature lines (cf. examples in
Section \ref{se.baeck})

\section{The conformal compactification of Euclidean space}

Let us begin by recalling some basic facts. The one-point
conformal compactification $\mathcal{M}=\R^{3}\cup\{\infty\}$ of
Euclidean space is classically realized as the projectivization of
the light cone $\mathcal{L}$ of Minkowski 5-space $\R^5_1$ with Lorentz
scalar product $\langle~,~\rangle$: $$ \mathcal{M}\cong \P(\mathcal{L})=
\left\{ [X]\in \P^4 : \langle X, X\rangle=0 \right\}. $$ We
consider linear coordinates $x^0,\dots,x^4$ such that
\begin{equation}
\langle X,Y\rangle=-\left(x^{0}y^{4}+x^{4}y^{0}\right)+x^{1}y^{1}+
x^{2}y^{2}+ x^{3}y^{3}=g_{ij}x^{i}y^{j}, \label{innerproduct}
\end{equation}
and identify $\mathcal{M}$ with $\P [\mathcal{L}]$ by means of the
conformal map $$ j:(p^{1},p^{2},p^{3})\in \R^{3}\mapsto
\left[\left(1,p^{1},p^{2},p^{3},\frac{1}{2}|p|^{2}\right)\right]\in
\P [\mathcal{L}], \quad
j(\infty)=\left[\left(0,0,0,0,1\right)\right]. $$ The linear
action of the pseudo-orthogonal group $G\cong SO(4,1)$ of
(\ref{innerproduct}) descends to a transitive action on $\mathcal{M}$
by conformal (M\"obius) diffeomorphisms. In this model for $\mathcal{M}$, the
de~Sitter space $S^4_1=\left\{Y\in \R^5_1 : \langle Y,Y\rangle =
1\right\}$ parametrizes the 2-spheres in $\mathcal{M}$ by $$ Y\mapsto
\P((Y)^{\perp}\cap\mathcal{L}) $$ and $G$ acts transitively on the set
of 2-spheres.

\medskip
\noindent A {\it M\"obius frame} is a basis $B=(B_0,\dots,B_4)$ of
$\R_1^5$ such that the vectors form the columns of a matrix of
$G$. Geometrically, the unit space-like vectors $B_1$, $B_2$,
$B_3$ represent 2-spheres which intersect orthogonally, and
$[B_0],[B_4]\in \mathcal{M}$ their intersection points. Regarding
$B_0,\dots,B_4$ as $\R^5$-valued functions defined on $G$, there
are unique 1-forms $\{\omega^{I}_{J}\}_{0\leq I,J\leq 4}$ such
that $$ dB_I=\omega^{J}_{I}B_J,\quad
\omega^{K}_{I}g_{IJ}+\omega^{K}_{J}g_{KI}=0 $$ and satisfying the
structure equations $$ d\omega^{I}_{J} = - \omega^{I}_{K} \wedge
\omega^{K}_{J}. $$ $(\omega^{I}_{J})=B^{-1}dB$ is the
Maurer-Cartan form of $G$ with values in the Lie algebra
$\mathcal{G}$ of $G$.

\begin{remark} 
The set $\dot\mathcal{L}$ of all $X\in \mathcal{L}$ such that
$x^{4}\neq 0$ can be given a Lie group structure. For, let $X\in
\dot\mathcal{L}$ and define $$ g^{+}(X)=\left(\begin{array}{ccccc}
 1/x^{4}&  x^{1}/x^{4} & x^{2}/x^{4}& x^{3}/x^{4} & x^{0}\\
         0 & 1 & 0 & 0 & x^{1} \\
          0 & 0 & 1 & 0 & x^{2}\\
        0 & 0 & 0 & 1 & x^{3} \\
         0 & 0 & 0 & 0 & x^{4}
\end{array} \right).
$$ The mapping $$ g^{+}:X\in\dot{\mathcal{L}}\mapsto g^{+}(X)\in G $$
is a smooth embedding which induces a Lie group structure on
$\dot\mathcal{L}$. The group operation is given by $$ X\star Y=\left(
x^{0}y^{4}+\frac{1}{x^{4}}(y^{0}+\sum_{j=1}^{3}x^{j}y^{j}),
x^{1}+x^{4}y^{1},
x^{2}+x^{4}y^{2},x^{3}+x^{4}y^{3},x^{4}y^{4}\right)^t. $$ In
particular, $X^{-1}=(x^{0},\frac{-x^1}{x^{4}},
\frac{-x^2}{x^{4}},\frac{-x^3}{x^{4}},\frac{1}{x^{4}}{)^t}$ and
$1=(0,0,0,0,1{)^t}$. 
\end{remark}

\section{The Lax pair} \label{se.lax}

\begin{defn}
Throughout the paper a solution $\Phi$ of equation
(\ref{solitoneq}) will be referred to as a {\it wave potential}
with {\it character} $s$.
\end{defn}

\noindent Let $U\subset \R^2$ be a simply connected domain with
coordinates $(x,y)$, and let $\Phi :U\to \R$ be a nowhere
vanishing smooth function. For $m\in \R$, we define $P(m),Q(m) :U
\to \mathcal{G}$ by
\begin{eqnarray*}
 P(m)&=&\left(\begin{array}{ccccc}
          2\Phi^{-1}\Phi_{x}& m\Phi^{-1}-\frac{s}{2}\Phi^{-3}&
           0&   0 & 0\\
         \Phi & 0 &-\Phi^{-1}\Phi_{y}& -\Phi&
         m\Phi^{-1}-\frac{s}{2}\Phi^{-3} \\
          0 & \Phi^{-1}\Phi_{y}& 0 & 0 & 0\\
 0 & \Phi& 0  & 0 & 0   \\
 0 & \Phi & 0 & 0 & -2\Phi^{-1}\Phi_{x}
\end{array} \right),
\\
Q(m)&=&\left(\begin{array}{ccccc}
          2\Phi^{-1}\Phi_{y}& 0 &
           -m\Phi^{-1}-\frac{s}{2}\Phi^{-3}& 0 & 0\\
         0 & 0 & \Phi^{-1}\Phi_{x} & 0 & 0 \\
          \Phi& -\Phi^{-1}\Phi_{x} & 0 & \Phi &
          -m\Phi^{-1}-\frac{s}{2}\Phi^{-3}\\
 0 & 0 & -\Phi & 0 & 0   \\
 0 & 0 & \Phi& 0 & -2\Phi^{-1}\Phi_{y}
\end{array} \right).
\end{eqnarray*}

A straightforward calculation gives:

\begin{lemma}
Equation (\ref{solitoneq}) is equivalent to the matrix Lax
equation
\begin{equation}
 P(m)_y - Q(m)_x - [P(m),Q(m)] = 0,
\end{equation}
where $P(m), Q(m)$ are as defined above.
\end{lemma}

\noindent This is the compatibility condition for the linear
system
\begin{equation}
V_x=-P(m)V, \quad V_y= -Q(m)V\label{darbouxsys}
\end{equation}
in Minkowski 5-space. The linear system (\ref{darbouxsys}) is
referred to as the $D_m$-{\it system}\footnote{The linear system
(\ref{darbouxsys}) is gauge equivalent to a special case of the
system considered by Darboux and Bianchi for the construction the
Darboux transformation \cite{B1, BHPP} (see also Section
\ref{se.baeck})} associated to $\Phi$. The $\mathcal{G}$-valued
one-form $\beta(m):=P(m)dx+Q(m)dy$ satisfies the Maurer-Cartan
equation $d\beta(m) +\beta(m)\wedge\beta(m)=0$ and can
consequently be integrated to a frame $B(m) : U\to G$, uniquely
determined up to left multiplication by a constant element of $G$,
such that\
\begin{equation}
dB(m)=\beta(m)B.\label{structureequation}
\end{equation}
We call $B(m)$ a {\it normal frame field} for the wave potential
$\Phi$ with {\it spectral parameter} $m$. 
Because of (\ref{structureequation}), the components $v^{0},\dots,v^{4}$ of 
any constant vector $X\in
\R^{5}_{1}$ with respect to a normal frame $B(m)$ provide a
solution $V$ of the $D_{m}$-system. From (\ref{structureequation})
it also follows that
\begin{equation}
b(m)= \frac{s}{2}\Phi^{-2}B(m)_{0}+mB(m)_{3}+\Phi^{2}B(m)_{4}
\end{equation}
is a constant vector; we call $b(m)$ the {\it pointing vector} of
the normal frame $B(m)$. The {\it kinetic energy} $C_{1}(V)$ and
the {\it linear momentum} $C_{2}(V)$ of a solution $V$ of the
linear system (\ref{darbouxsys}) are defined by
\begin{eqnarray}
C_{1}(V) &=&-2v^{0}v^{4}+(v^{1})^{2}
+(v^{2})^{2}+(v^{3})^{2},\label{conserv1}
\\
C_{2}(V)&=&\Phi^{2}v^{0}-mv^{3}+
\frac{s}{2}\Phi^{-2}v^{4}.\label{conserv2}
\end{eqnarray}
If $X\in \R^{5}_{1}$ is the initial condition of $V$ with respect
to the normal frame $B(m)$, then $C_{1}(V)= \langle X,X \rangle$
and $C_{2}(V)=-\langle b(m),X\rangle$. Thus, $C_{1}$ and $C_{2}$
are two first integrals of the $D_{m}$-system.

\section{Potentials and special isothermic surfaces}\label{se.geom}

\noindent Let $f : U\to \R^{3}$ be an isothermic immersion with
conformal principal coordinates $x,y$ and let
\begin{equation}
I=\Theta^{2}(dx^{2}+dy^{2}),\quad
II=\Theta^{2}(h_{1}dx^{2}+h_{2}dy^{2}), \label{Fundamentalforms}
\end{equation}
be its fundamental forms, where $\Theta$ is a nowhere vanishing
smooth function and $h_{1}$ and $h_{2}$ are the principal
curvatures.

\begin{remark}
Recall that the notion of an isothermic immersion is
conformally invariant, that is, if $f$ is an isothermic immersion
and $A\in G$ is a conformal diffeomorphism, then $A\circ f$ is
also isothermic \cite{Bu,M2}. In the following we will not make any
distinction between isothermic immersions in $\R^{3}$ or in 
$\mathcal{M}$. 
\end{remark}


\noindent Following Bianchi \cite{B1} we set $$
2H=h_{1}+h_{2},\quad L=\Theta^{2}(h_{1}-h_{2}),\quad M=-HL. $$ and
give the following
\begin{defn}
$f$ is called {\it special isothermic} of type $(A,B,C,D)$ if
there exist real constants $A,B,C,D$ such that
\begin{equation}
4\Theta^{2}|\nabla H|^{2}+M^{2}+2AM +2BH
+2CL+D=0.\label{specialeq}
\end{equation}
\end{defn}

\noindent We shall assume that $f:U\to \R^{3}$ is {\it generic},
i.e., either one of the following conditions hold: $L_{x}L_{y}\neq
0$, or $L_{x}\neq 0$ and $ L_{y}= 0$, or $L_{y}\neq 0$ and
$L_{x}=0$, or else $L_{x}=L_{y}=0$.

\medskip
\noindent Next, let
$
\Phi=\frac{1}{2}\left(h_{1}-h_{2}\right)\Theta
$
be the {\it Calapso potential} \footnote{For more information on
the Calapso potential see \cite{C1,Be,Bu,BHPP,M2}} of the
isothermic immersion $f$.

\medskip
We are now in a position to state:

\begin{prop}
Let $f:U\to \R^{3}$ be a special isothermic immersion of type
$(A,0,C,D)$. Then the Calapso potential $\Phi$ is a wave potential
with character $D/4$.
\end{prop}

\noindent \begin{proof}{\rm The Gauss and Codazzi equations
\begin{equation}
\Theta\Delta\Theta-|\nabla
\Theta|^{2}+h_{1}h_{2}\Theta^{4}=0,\quad \Theta\neq 0,
\label{Gaussequation}
\end{equation}
\begin{equation}
(h_{1})_{y}=-\frac{\Theta_{y}}{\Theta}(h_{1}-h_{2}),\quad
(h_{2})_{x}=\frac{\Theta_{x}}{\Theta}(h_{1}-h_{2}),
\label{Mainardiequation}
\end{equation}
imply 
$$ 
L_{x}=2\Theta^{2}H_{x},\quad
L_{y}=-2\Theta^{2}H_{y},\quad
L_{xy}=\Theta^{-1}\left(\Theta_{y}L_{x}+\Theta_{x}L_{y}\right),
$$
and 
$$ 
M_{x}=-h_{1}L_{x},\quad M_{y}=-h_{2}L_{y}. 
$$ 
Then, (\ref{specialeq}) becomes \footnote{Note that the functions
$\Theta^{-1},-h_{1}\Theta^{2},h_{2}\Theta^{2}$ are related to the
Christoffel transformation of $f$. Thus, from (\ref{varspecialeq})
it follows that the Christoffel transform of an isothermic
immersion of type $(A,0,C,D)$ is special of type $(A,C,0,D)$, see
also \cite{B1}.}
\begin{equation}
\Theta^{-2}|\nabla L|^{2}+M^{2}+2AM+2CL+D=0.\label{varspecialeq}
\end{equation}
Since
\begin{equation}
2\Phi = (h_{1}-h_{2})\Theta,\quad L=2\Phi \Theta,\quad M=-2\Phi H
\Theta,\label{varhlm}
\end{equation}
equation (\ref{varspecialeq}) also reads
\begin{eqnarray}
\Theta^{-2}|\nabla \Theta |^{2}&=&-\Phi^{-2}|\nabla
\Phi|^{2}-2\Phi^{-1}\Theta^{-1}\nabla \Phi\cdot \nabla \Theta
+A\Phi^{-1}H\Theta  \nonumber\\ & & {}-H^{2}\Theta^{2}
-C\Phi^{-1}\Theta-\frac{D}{4}\Phi^{-2}. \label{gradtheta}
\end{eqnarray}
The derivative of (\ref{varspecialeq}) with respect to $x$ and $y$
yields
\begin{eqnarray}
2\Theta^{-1}L_{x}\left((\Theta^{-1}L_{x})_{x}+
\Theta^{-2}L_{y}\Theta_{y}-\Theta(M+A)h_{1}+
C\Theta\right)&=&0,\nonumber\\
2\Theta^{-1}L_{y}\left((\Theta^{-1}L_{y})_{y}+
\Theta^{-2}L_{x}\Theta_{x}-\Theta(M+A)h_{2}+
C\Theta\right)&=&0.\nonumber
\end{eqnarray}
When $L_{x}L_{y}\neq 0$ the last two equations imply $$
(\Theta^{-1}L_{x})_{x}+(\Theta^{-1}L_{y})_{y}
+\Theta^{-2}(L_{x}\Theta_{x}+L_{y}\Theta_{y}) -2\Theta
H(M+A)+2C\Theta=0. $$ This equation combined with (\ref{varhlm})
and (\ref{gradtheta}) gives the result. A similar argument applies
in the other cases. }\end{proof}


\begin{remark}
In the setting of the previous proposition, notice that when
$L_x=L_y=0$ the immersion $f$ has constant mean curvature. When
instead $L_x=0$ and $L_y\neq 0$ the Calapso potential
$\Phi=\Phi(y)$ is a function of the variable $y$ alone and $f$ is
an isothermic canal surface. By a classical result of Darboux
\cite{D}, $f$ is then conformally equivalent to either a cone, a
cylinder, or a surface of revolution. If $f$ is a surface of
revolution, one can prove that the profile curve lies in the
hyperbolic half-plane, that its arclength is proportional to $y$,
and that its curvature is parametrized by the  Calapso potential.
Further, equation (\ref{solitoneq}) reduces to the Euler-Lagrange
equation of the total square curvature functional $\alpha
\mapsto\frac{1}{2}\int{\kappa^2_\alpha}$ on smooth curves $\alpha$
with fixed length. Therefore, the surface defined by $f$ is
obtained by revolving an elastic curve (possibly a free elastic
one) in the hyperbolic half-plane about the boundary at infinity.
This description indicates, in particular, how to construct
special isothermic surfaces with umbilic lines (see also
Babich--Bobenko \cite{BB}) Recall that rotational Willmore
isothermic surfaces arise from free elastic curves in the
hyperbolic 2-plane \cite{BG}.  See also Langer--Singer \cite{LS}
and Pinkall \cite{Pi}. For the case of Willmore  canal surfaces we
refer to \cite{MNts}. 
\end{remark}

As for the converse, consider first the following:

\begin{defn} 
For any $X\in \dot\mathcal{L}$, let 
$\pi_{X}:\mathcal{L}\to \mathcal{M}$ 
be the conformal map defined by $\pi_{X}=j^{-1}\circ \pi\circ
g^{+}(X)^{-1}$, where $\pi : \mathcal{L}\to \P [\mathcal{L}]$ denote the
canonical projection. 
\end{defn}

\begin{thm}
Let $B(m) : U\to G$ be a normal frame field with wave potential
$\Phi$ and character $s$, and let $X\in \dot\mathcal{L}$ be a
constant vector. Then, $f_m=\pi_{X}\circ B(m)_{0}:U\to \mathcal{M}$ is a
special isothermic immersion of type $(-2m,0,2\langle
b(m),X\rangle,4s)$ with Calapso potential $\Phi$.
\end{thm}

\noindent \begin{proof}{\rm Let $V$ be the solution of the
$D_{m}$-system with initial condition $X$, that is $X=B({m})V$.
Then, the fundamental forms of $f_m =\pi_{X}\circ B(m)_{0}$ can be
read off the Maurer--Cartan form of the Euclidean frame
$B(m)g^{+}(V)$. These are computed to be
\begin{equation}
I =
\frac{\Phi^{2}}{(v^{4})^{2}}\left((dx)^{2}+(dy)^{2}\right),
\quad
II  =
\frac{\Phi^{2}}{(v^{4})^{2}}\left((v^{4}-v^{3})(dx)^{2}-(v^{4}+v^{3})(dy)^{2}\right).
\end{equation}
This implies that $f_m$ is an isothermic immersion and that
$$ 
H = -v^{3},\quad L = \frac{2}{v^{4}}\Phi^{2},\quad M =
2\frac{v^{3}}{v^{4}}\Phi^{2}. 
$$ 
Next, by using the constraint
$\langle V,V\rangle=0$ and the conservation of the linear momentum
$C_{2}(V)=-\langle b(m),X\rangle$, a direct computation shows that
\begin{equation}
4\Theta^{2}|\nabla H|^{2} + M^{2} - 4mM + 4\langle b(m),X\rangle L
+ 4s=0.
\end{equation}
Thus $f$ is special isothermic of type $(-2m,0,2\langle
b(m),X\rangle,4s)$.}
\end{proof}

\begin{remark}
The solution $V$ in the proof of the theorem can be expressed
in terms of the Euclidean invariants of $f$ by the following
formulae:
\begin{eqnarray*}
v^{0}
=\frac{1}{h_{1}-h_{2}}\left(\frac{4(H_{x})^{2}}{(h_{1}-h_{2})^{2}\Theta^{2}}+
\frac{4(H_{y})^{2}}{(h_{1}-h_{2})^{2}\Theta^{2}}+H^{2}\right),\\
v^{1}=\frac{2H_{x}}{(h_{1}-h_{2})\Theta},\quad
v^{2}=-\frac{2H_{y}}{(h_{1}-h_{2})\Theta}, \quad v^{3}=-H, \quad
v^{4}=\frac{1}{2}(h_{1}-h_{2}).
\end{eqnarray*}
\end{remark}

\begin{remark}
Given an umbilic free immersion $f:U\to \R^{3}$, there exist
a canonical frame field $B_{f}:U\to G$ along $f$: the {\it central
frame} field of the immersion. The construction of such a frame is
due to Bryant \cite{Br}. If $f$ is special isothermic of type
$(A,0,C,D)$ with Calapso potential $\Phi$, then $B_{f}$ is exactly
the normal frame field with wave potential $\Phi$ and spectral
parameter $m=-A/2$. 
\end{remark}

\section{The B\"acklund transformation and $1$-soliton solutions}

\begin{defn}
Let $\Phi : U\to \R$ be a wave potential with character
$s$. A solution $V$ of the linear system (\ref{darbouxsys})
satisfying $C_{1}(V)=C_{2}(V)=0$ is said to be an {\it $m$-system
of transforming functions} \footnote{Note that $m$-transforming
functions do exist only if $m^{2}-s\geq 0$. If $m^{2}-s>0$, then
the set of $m$-transforming functions with potential $\Phi$ is a
$3$-dimensional cone. If $m^{2}-s =0$, the transforming functions
are of the form 
$$
V=r\left(\frac{s}{2}\Phi^{-2},0,0,m,\Phi^{2}\right)^t 
$$ 
for a constant $r\neq 0$. Observe that the last component $v^{4}$ of a
system of transforming functions never vanishes. If, in addition,
the character $s$ is positive, then also the third component
$v^{3}$ never vanishes.} for the potential $\Phi$. The
corresponding {\it B\"acklund transform} is defined by

\begin{equation}
E(\Phi,V)=\frac{v^{3}}{v^{4}}\Phi \label{Backlundtransform}
\end{equation}
\end{defn}

\begin{figure}[ht]
\begin{center}
\fbox{
\includegraphics[height=6cm,width=9cm]{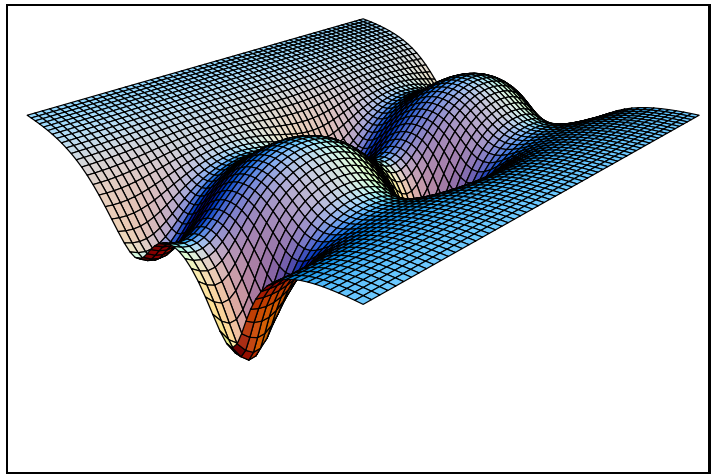}}
\caption{a B\"acklund transform $\Psi(x,y)$ of the wave potential
$\sqrt{2}Sech(\sqrt{2}y)$.}
\end{center}
\end{figure}

\begin{thm}
The B\"acklund transform $E(\Phi,V)$ is a wave potential with
character $s$.
\end{thm}

\noindent\begin{proof}{\rm Let $\tilde{\Phi}=E(\Phi,V)$. The first
derivatives of $\tilde{\Phi}$ are given by $$
\tilde{\Phi}_{x}=\frac{v^3}{v^4}\left[-\Phi_{x}+v^{1}\left(
\frac{1}{v^{4}}-\frac{1}{v^{3}}\right){\Phi}^2\right],\quad
\tilde{\Phi}_{y}=\frac{v^3}{v^4}\left[-\Phi_{y}+v^{2}\left(
\frac{1}{v^{4}}+\frac{1}{v^{3}}\right){\Phi}^2\right], $$ and then
\begin{eqnarray*}
\tilde{\Phi}_{xx} & = &\frac{v^3}{v^4}\left\{-\Phi_{xx}
+2{\Phi_x^2}{\Phi}^{-1}
+\Phi^2\left(\frac{1}{v^{4}}-\frac{1}{v^{3}}\right)
\left[-{\Phi}v^{0}-\Phi_{x}{\Phi}^{-1}v^{1}
+\Phi_{y}{\Phi}^{-1}v^{2} \right.\right.
\\
& & \left.\left. {} + {\Phi}v^{3} + \left(\frac{s}{2}{\Phi}^{-3}
-m{\Phi}^{-1}\right)v^{4} \right] -2{\Phi}{\Phi_x}\frac{v^1}{v^4}
-2{\Phi}^{3}\frac{({v^1})^2}{v^{3}v^{4}}+
2{\Phi}^{3}\left(\frac{v^{1}}{v^{4}}\right)^{2}\right\},
\\
\tilde{\Phi}_{yy} & = &\frac{v^3}{v^4}\left\{-\Phi_{yy}
+2{\Phi_y^2}{\Phi}^{-1}
+\Phi^2\left(\frac{1}{v^{4}}+\frac{1}{v^{3}}\right)
\left[-{\Phi}v^{0}+\Phi_{x}{\Phi}^{-1}v^{1}
-\Phi_{y}{\Phi}^{-1}v^{2} \right.\right.
\\
& & \left.\left. {} + {\Phi}v^{3} + \left(\frac{s}{2}{\Phi}^{-3}
-m{\Phi}^{-1}\right)v^{4} \right] -2{\Phi}{\Phi_y}\frac{v^2}{v^4}
-2{\Phi}^{3}\frac{({v^2})^2}{v^{3}v^{4}}+
2{\Phi}^{3}\left(\frac{v^{2}}{v^{4}}\right)^{2}\right\}.
\end{eqnarray*}
We now use the constraint $\langle V,V\rangle =0$ to obtain $$
{\tilde{\Phi}}^{-1}\Delta\tilde{\Phi}-{\tilde{\Phi}}^{-2}|\nabla{\tilde{\Phi}}|^2
=
-{\tilde{\Phi}}^{2}+\frac{v^{4}}{(v^{3})^{2}}
\left(2mv^{3}-2{{\Phi}}^{2}v^{0}\right). $$ From this equation and
the constraint $$
{\Phi}^{2}v^{0}-mv^{3}+\frac{s}{2}{\Phi}^{-2}v^{4}=0 $$ the result
follows. }
\end{proof}

\begin{remark} 
Let $\Phi$ be a potential with character $s\neq
0$, then the {\it complementary potential} $\Phi^{\ast}$ is
defined by $\sqrt{|s|}\Phi^{-1}$. This is a new potential with the
same character of $\Phi$. Note that $\Phi^{\ast}$ can be obtained
as the B\"acklund transform of $\Phi$ with respect to the system
of transforming functions 
$$
V=\left(\frac{s}{2}\Phi^{-2},0,0,\sqrt{|s|},\Phi^{2}\right)^t. 
$$
\end{remark}

 \begin{ex}[\textbf{One-soliton solutions}] 
Consider the trivial
 solution $\Phi = 1$ of (\ref{solitoneq}) with  $s=1$.
In this case the $1$-form
 $\beta({m})$ corresponding to the Lax pair is given by

 \begin{equation}
 \beta({m})=\left(\begin{array}{ccccc}
           0& (m-\frac{1}{2})dx & -(m+\frac{1}{2})dy &   0 & 0\\
          dx& 0 & 0 & -dx &(m-\frac{1}{2})dx \\
           dy& 0 &0& dy &-(m+\frac{1}{2})dy\\
           0& dx& -dy & 0&0   \\
  0&dx&dy&0&0
 \end{array} \right).
 \end{equation}
Now set $$ \zeta = \sqrt{2| m-1 |},\quad \eta = \sqrt{2 |m+1|}. $$
By solving a system of first order linear differential equations
with constant coefficients, the normal framing $B(m)$ is computed
to be $$ B(m)=\left(\begin{array}{ccccc} \frac{\eta \cosh(\zeta x
 )+2}{\sqrt{2}\zeta\eta} & \frac{\sinh(\zeta x)}{\sqrt{2}}& 0 &
 -\frac{\eta\cosh(\zeta x)+2m}{\sqrt{2}\eta\zeta}
 &\frac{(2m-1)\eta\cosh(\zeta x)+2}{2\sqrt{2}\eta\zeta}\\
         \frac{\sinh(\zeta x)}{\zeta} & \cosh(\zeta
 x) & 0 & -\frac{\sinh(\zeta x)}{\zeta} &\frac{(2m-1)\sinh(\zeta x)
 }{2\zeta}\\
          \frac{\cos(\eta y)}{\eta} & 0 &-\sin(\eta y)& \frac{\cos(\eta y)}{\eta}
          &-\frac{(2m+1)\cos(\eta y)}{2\eta}\\
          \frac{\sin(\eta y) }{\eta} & 0& \cos(\eta y) & \frac{\sin(\eta y)}{\eta}&
          -\frac{(2m+1)\sin(\eta y)}{2\eta}   \\
  \frac{\eta \cosh(\zeta
 x)-2}{\sqrt{2}\zeta\eta}&\frac{\sinh(\zeta
 x)}{\sqrt{2}}&0&-\frac{\eta\cosh(\zeta
 x)-2m}{\sqrt{2}\eta\zeta}&\frac{(2m-1)\eta\cosh(\zeta
 x)-2}{2\sqrt{2}\eta\zeta}
 \end{array} \right),
$$ 
In order to find the solutions of (\ref{darbouxsys}) satisfying
the constraint $C_{1}=C_{2}=0$, we may assume that $m>1$. Note
that the pointing vector of the normal framing $B(m)$ is the
space-like vector $(1,0,0,0,-1)^t$. Thus, the transforming functions
of the $D_{m}$-system are given by 
$$
V(m,a^{1},a^{2},a^{3})=G\cdot B({m})^t\cdot GX(a^{1},a^{2},a^{3}),
$$ where $X(a^{1},a^{2},a^{3})\in \mathcal{L}$ is defined by 
$$
{X}(a^{1},a^{2},a^{3})=\left(
\sqrt{\frac{(a^{1})^{2}+(a^{2})^{2}+(a^{3})^{2}}{2}}
,a^{1},a^{2},a^{3},\sqrt{\frac{(a^{1})^{2}+(a^{2})^{2}+(a^{3})^{2}}{2}}
\right)^t 
$$ and where $G=(g_{ij})$ is as in (\ref{innerproduct}).
It follows that
\begin{eqnarray}
v^{0}&=&\frac{2m+1}{2\eta}\left( a^{2}
\cos(\eta y)+a^{3} \sin (\eta
y)\right)+\frac{2m-1}{\sqrt{2}\eta\zeta}(\|a\|\cosh(\zeta
x)+\nonumber\\
& & {}-\frac{a^{1}\eta}{\sqrt{2}}\sinh(\zeta x))\nonumber\\
v^{1}&=&a^{1}\cosh(\zeta x)-\|a\|\sinh(\zeta x)\nonumber\\
v^{2}&=&a^{1}\cos(\eta y)-a^{2}\sin(\eta y)\\
v^{3}&=&\frac{a^{2}}{\eta}\cos(\eta y)+\frac{a^{3}}{\eta}\sin(\eta
y)+\frac{\|a\|}{\zeta}\cosh(\zeta
x)-\frac{a^{1}}{\zeta}\sinh(\zeta x)\nonumber\\
v^{4}&=&\frac{\|a\|}{\zeta}\cosh(\zeta
x)-\frac{a^{1}}{\zeta}\sinh(\zeta x)-\frac{a^{2}}{\eta}\cos(\eta
y)-\frac{a^{3}}{\eta}\sin(\eta y),\nonumber
\label{transformingfunctions}
\end{eqnarray}
where 
$$
\parallel
a\parallel=\sqrt{(a^{1})^{2}+(a^{2})^{2}+(a^{3})^{2}}. 
$$ 
From this we obtain the following formula for the one-soliton solutions
\begin{equation}
\begin{split} 
&\Phi(m,a^{1},a^{2},a^{3})=\\
&\frac{\zeta(a^{2}\cos(\eta
y)+a^{3}{\eta}\sin(\eta y))+\eta(\|a\|\cosh(\zeta
x)-a^{1}\sinh(\zeta x))}{-\zeta(a^{2}\cos(\eta
y)+a^{3}{\eta}\sin(\eta y))+\eta( \|a\|\cosh(\zeta
x)-a^{1}\sinh(\zeta x))}.
\end{split}
\label{onesoliton}
\end{equation}

\begin{figure}[ht]
\begin{center}
\fbox{
\includegraphics[height=6cm,width=7cm]{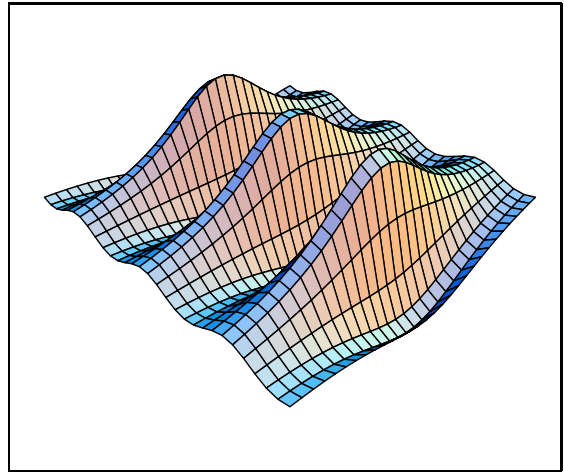}}
\caption{
one-soliton $\Phi(10/9,0,1/500,0)$.}
\end{center}
\end{figure}

\end{ex}

\section{The geometry of the B\"acklund transformation}\label{se.baeck}

\noindent In this section we briefly describe the geometric
transformation of special isothermic surfaces corresponding to the
B\"acklund transformation.


\begin{defn} 
A {\it curved flat
\footnote{For the general notion of a curved flat we refer to
Ferus--Pedit \cite{FP}.} framing} is a smooth map $A:U\to G$ such
that $A^{-1}dA$ takes the form $$ \left(\begin{array}{ccccc}
          0& -r\Theta^{-1}dx &
           r\Theta^{-1}dy&   0 & 0\\
\Theta dx& 0 &\Theta^{-1}(\Theta_{y}dx-\Theta_{x}dy) &
-h_{1}\Theta dx &
         -r\Theta^{-1}dx \\
   \Theta dy& \Theta^{-1} (\Theta_{x}dy-\Theta_{y}dx) & 0 &-h_{2}\Theta dy&
          r\Theta^{-1}dy\\
 0& h_{1}\Theta dx& h_{2}\Theta dy & 0&0   \\
 0&\Theta dx&\Theta dy&0&0
\end{array} \right),
$$ where $\Theta$, $h_{1}$, $h_{2}$ are smooth functions with
$\Theta(p)\neq 0$, for each $p\in U$, and $r$ is a constant
referred to as the {\it spectral parameter} of the curved flat.
For short, the connection form $A^{-1}dA$ will be denoted by
$\alpha_{r}(\Theta,h_{1},h_{2})$. 
\end{defn}

\noindent The functions $\Theta$, $h_{1}$ and $h_{2}$ satisfy the
{\it isothermic Gauss-Codazzi system}
\begin{equation}
\left\{ \begin{array}{lll} 
& \Theta\Delta(\Theta)=|\nabla
\Theta|^{2}-h_{1}h_{2}\Theta^{4},\\
& \Theta (h_{1})_{y}=\Theta_{y}(h_{2}-h_{1}),\\
& \Theta (h_{2})_{x}=\Theta_{x}(h_{1}-h_{2}).\\
\end{array}\right. \label{GMsystem}
\end{equation}

\begin{figure}[ht]
\begin{center}
\fbox{
\includegraphics[height=6.5cm,width=6.5cm]{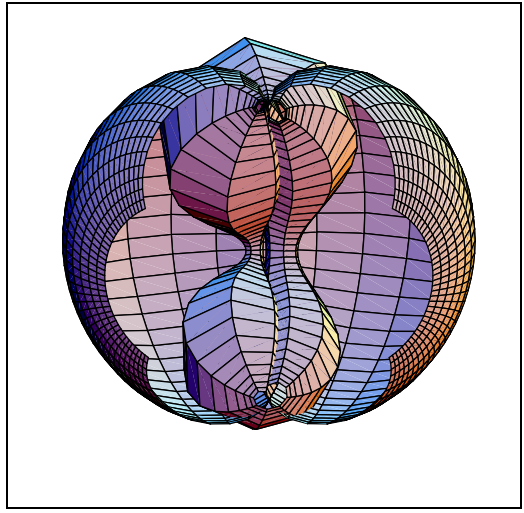}}
\caption{special isothermic surface with wave potential
$\Psi(x,y)$ (see Figure $1$) and spectral parameter $m=1$.}
\end{center}
\end{figure}


\begin{remark}
Note that if $A:U\to G$ is a curved flat framing, then
 $f=[A_{0}]:U\to \mathcal{M}$ and $\bar{f}=[A_{4}]:U\to
\mathcal{M}$ are isothermic immersions with Calapso potentials $$
\Phi=\frac{1}{2}(h_{1}-h_{2})\Theta,\quad
\bar{\Phi}=-\frac{1}{2}(h_{1}+h_{2})\Theta, $$ respectively. We
say that $A$ is a curved flat framing along $f$ and call $\Phi$
the Calapso potential of $A$. The {\it dual frame field} is
defined to be $\bar{A}=AJ$, where $J\in G$ is given by
\begin{equation}
\left(\begin{array}{ccccc}
 0& 0& 0& 0 & 1\\
 0 & -1 &0& 0& 0 \\
 0 & 0& 1 &0& 0\\
 0& 0& 0 & 1&0   \\
 1 & 0 & 0 & 0 & 0
\end{array} \right).
\end{equation}
In particular, $\bar{A}$ is a curved flat framing along $\bar{f}$.
The mapping $\sigma:=A_{3}:U\to S^{4}_{1}\subset \R^{5}_{1}$
defines a sphere congruence whose envelopes are $f$ and $\bar{f}$.
Moreover, the correspondence induced by $\sigma$ preserves the
curvature lines and is conformal. This amounts to saying that
$\sigma$ is a Darboux congruence. Accordingly, $\bar{f}$ is called
a {\it Darboux transform} of $f$. Any Darboux transform of $f$
arises in this way \footnote{The result that an isothermic surface
together with a Darboux transform form a curved flat in the
pseudo-Riemannian symmetric space of pairs of distinct points in
$\mathcal{M}$ has been proved in \cite{BHPP}. For the interpretation of
the Darboux transforms as dressing transformations of loop groups
we refer to Burstall \cite{Bu}}. 
\end{remark}

\begin{remark}
The first envelope $f:U\to \mathcal{M}$ of a
curved flat framing $A$ is a special isothermic immersion if and
only if, in addition to (\ref {GMsystem}), the triple
$(\Theta,h_{1},h_{2})$ satisfies an equation of the form
\begin{equation}
4\Theta^{2}|\nabla H|^{2}+M^{2}+2AM+2CL+D=0,
\label{specialcurvedflatequation}
\end{equation}
 where $A$, $C$, $D$ are real constants and $M=-HL$ being
$L=\Theta^{2}(h_{1}-h_{2})$ and $2H=h_{1}+h_{2}$. We then say that
$A$ is a {\it special curved flat framing} of type $(A,0,C,D)$.
Notice that both the envelopes of the curved flat are special
isothermic immersions if and only if $C=0$. In this case we shall
say that $\bar{f}$ is a {\it special Darboux transform} of $f$. 
\end{remark}

We now have:

\begin{prop}
Let $f:U\to \mathcal{M}$ be a special isothermic immersion of type
$(A,0,C,D)$ with deformation parameter $m$ and Calapso potential
$\Phi$, and let $B(m):U\to G$ be the Bryant's central frame field
along $f$. Next, let $V:U\to \mathcal{L}$ be any solution of the
$D_{h}$-linear system $$ dV=-\beta(h)V,\quad h\neq m, $$
satisfying the constraint $C_{2}(V)=0$. Then
 $$
\bar{f}:=\left[v^{0}B(m)_{0}+v^{1}B(m)_{1}+v^{2}B(m)_{2}+v^{3}B(m)_{
3}+v^{4}B(m)_{4}\right] $$ defines a special isothermic immersion
which is a special Darboux transform of $f$ such that $$
\bar{\Phi}=E(\Phi,V). 
$$
\end{prop}

\noindent \begin{proof}{\rm The connection form $\alpha$ of the
framing $A=B(m)g^{+}(V)$ $$
\alpha=g^{+}(V)^{-1}d[g^{+}(V)]+g^{+}(V)^{-1}\beta(m)g^{+}(V) $$
takes the form
$
\alpha=\alpha_{r}(\Theta,h_{1},h_{2}) $,
  with
\begin{equation}
\Theta=\frac{\Phi}{v^{4}},\quad h_{1}=v^{4}-v^{3},\quad
h_{2}=-(v^{4}+v^{3}),\quad r=h-m.
\end{equation}
 It then follows that $A$ is a special curved flat framing of
type $(-2h,0,-2C_{2}(V),4s)$, where $s$ denotes the character of
$\Phi$. The first envelope $[A_{0}]$ represents the original
immersion $f$ and the second envelope $[A_{4}]$ represents the
Darboux transform $\bar{f}$. If $C_{2}(V)=0$, then $\bar{f}$ is a
special isothermic immersion and $\bar{\Phi}=E(\Phi,V)$.
}\end{proof}

\begin{remark}
$\bar{f}$ has deformation parameter $m$ and
the normal frame field $\bar{B}(m):U\to G$ along $\bar{f}$ is
given by 
$$ 
\bar{B}(m)=B(m)g^{+}(V)Jg^{+}(T), 
$$ 
where $T:U\to \mathcal{L}$ 
is the smooth map defined by
\begin{equation}
T=\left(\frac{\Phi^{2}v^{0}v^{4}}{(h-m)(v^{3})^{3}},\frac{v^{1}v^{4}}{(v^{3})^{2}},
-\frac{v^{2}v^{4}}{(v^{3})^{2}},-\frac{v^{4}}{v^{3}},\frac{(h-m)\Phi^{-2}(v^{4})^{2}}{v^{3}}
\right)^t. \label{Tequation}
\end{equation}
\end{remark}

\begin{ex}[\textbf{Special Darboux transforms of Dupin cyclides}]
Special isothermic maps with Calapso potential $\Phi=1$ are
 given (up to the action of the conformal group) by the following
 formulae (where $m$ is the deformation parameter of the family)


 \begin{itemize}

 \item if $m>1$ :

 \begin{equation}
 f_{m}(x,y)=\left(\frac{\sqrt{2}\eta\sinh(\zeta x)}{2+\eta\cosh(\zeta x)},
 \frac{\sqrt{2}\zeta\cos(\eta y)}{2+\eta\cosh(\zeta x)},
 \frac{\sqrt{2}\zeta\sin(\eta y)}{2+\eta\cosh(\zeta x)}\right),
 \label{cone}
\end{equation}

\item if $m=1$ :

 \begin{equation}
 f_{1}(x,y)=\left(\frac{8x}{4x^{2}+1},\frac{4\cos(2y)}{4x^{2}+1},\frac{4\sin(2y)}{4x^{2}+1}\right),
 \label{cylinder}
\end{equation}

\item if $0\le m <1$ :

 \begin{equation}
 f_{m}(x,y)=\left(\frac{\sqrt{2}\eta\cos(\zeta x)}{2+\eta\sin(\zeta x)},
 \frac{\sqrt{2}\zeta\cos(\eta y)}{2+\eta\sin(\zeta x)},
 \frac{\sqrt{2}\zeta\sin(\eta y)}{2+\eta\sin(\zeta x)}\right).
 \label{torus}
\end{equation}

\end{itemize}

\noindent The surfaces $S_{m}\subset \R^{3}$ parametrized by the
maps $f_{m}:\R^{2}\to \R^{3}$ are the Dupin cyclides. It is a
classical result that Dupin cyclides are conformally equivalent 
to either a circular cone, a circular cylinder, or a torus of 
revolution \cite{V,Bl}.

\noindent The normal frame field along $f_{m}$
 is computed to be :
 \begin{itemize}

 \item if $m>1$ :

  $$ B(m)=\left(\begin{array}{ccccc} \frac{\eta \cosh(\zeta x
 )+2}{\sqrt{2}\zeta\eta} & \frac{\sinh(\zeta x)}{\sqrt{2}}& 0 &
 -\frac{\eta\cosh(\zeta x)+2m}{\sqrt{2}\eta\zeta}
 &\frac{(2m-1)\eta\cosh(\zeta x)+2}{2\sqrt{2}\eta\zeta}\\
         \frac{\sinh(\zeta x)}{\zeta} & \cosh(\zeta
 x) & 0 & -\frac{\sinh(\zeta x)}{\zeta} &\frac{(2m-1)\sinh(\zeta x)
 }{2\zeta}
 \\
          \frac{\cos(\eta y)}{\eta} & 0 &-\sin(\eta y)& \frac{\cos(\eta y)}{\eta}
          &-\frac{(2m+1)\cos(\eta y)}{2\eta}\\
          \frac{\sin(\eta y) }{\eta} & 0& \cos(\eta y) & \frac{\sin(\eta y)}{\eta}&
          -\frac{(2m+1)\sin(\eta y)}{2\eta}   \\
  \frac{\eta \cosh(\zeta
 x)-2}{\sqrt{2}\zeta\eta}&\frac{\sinh(\zeta
 x)}{\sqrt{2}}&0&-\frac{\eta\cosh(\zeta
 x)-2m}{\sqrt{2}\eta\zeta}&\frac{(2m-1)\eta\cosh(\zeta
 x)-2}{2\sqrt{2}\eta\zeta}
 \end{array} \right);
 $$

 \item if $m=1$ :

  $$ B(1)=\left(\begin{array}{ccccc}
\frac{1}{8}(4x^{2}+1)&x&0&-\frac{1}{8}(4x^{2}-3)&\frac{1}{16}(4x^{2}+9)\\
x&1&0&-x&\frac{x}{2}\\
 \frac{1}{2}\cos (2y)&0&-\sin (2y)&\frac{1}{2}\cos
(2y)&-\frac{3}{4}\cos (2y)\\
 \frac{1}{2}\sin (2y)&0&\cos
(2y)&\frac{1}{2}\sin (2y)&-\frac{3}{4}\sin (2y)\\
1&0&0&-1&\frac{1}{2}
 \end{array} \right);
 $$

 \item if $0\le m <1$ :

$$
  B(m)=\left(\begin{array}{ccccc}
  \frac{2+\eta \sin (\zeta x)}{\sqrt{2}\zeta \eta}&\frac{\cos(\zeta
  x)}{\sqrt{2}}&0&-\frac{2m+\eta \sin(\zeta x)}{\sqrt{2}\zeta
  \eta}&\frac{2+(2m-1)\eta\sin(\zeta x)}{2\sqrt{2}\zeta \eta}\\
  \frac{\cos(\zeta x)}{\zeta}&-\sin(\zeta x)&0&-\frac{\cos(\zeta
  x)}{\zeta}&\frac{(2m-1)\cos(\zeta x)}{2\zeta}\\
  \frac{\cos(\eta x)}{\eta}&0&-\sin(\eta y)&\cos(\eta y)&-\frac{(2m+1)\cos(\eta
  y)}{2\eta}\\
  \frac{\sin(\eta y)}{\eta}&0&\cos(\eta y)&\sin(\eta y)&\frac{(2m+1)\sin(\eta
  y)}{2\eta}\\
  \frac{2-\eta\sin(\zeta x)}{\sqrt{2}\zeta \eta}&\frac{\cos(\zeta
  x)}{\sqrt{2}}&0&\frac{-2m+\eta\sin(\zeta x)}{\sqrt{2}\zeta\eta}&\frac{2-(2m-1)\eta\sin(\zeta x)}{2\sqrt{2}\zeta\eta}
 \end{array} \right).
$$
\end{itemize}
 
\begin{figure}[ht]
\begin{center}
\fbox{
\includegraphics[height=6.5cm,width=6.5cm]{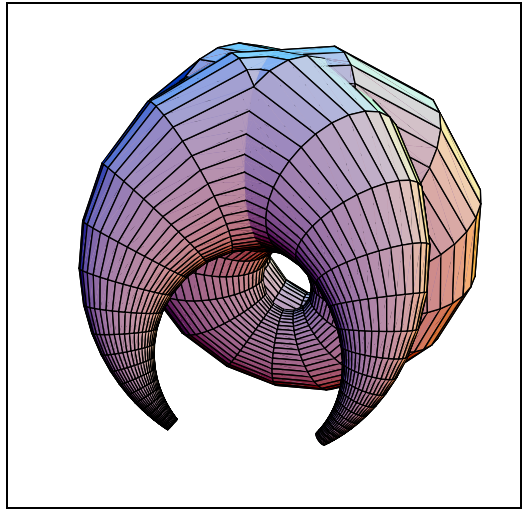}}
\caption{Darboux transforms of Dupin cyclides; special isothermic
surface with wave potential $\Phi(4/3,0,1/100,1/300)$ and spectral
parameter $m=1$.}
\end{center}
 \end{figure}


\noindent According to the above discussion, the special Darboux
transforms of $f_{m}$ are given by
\begin{equation}
 D(h,a^{1},a^{2},a^{3})(f_{m})=\left[B(m)V(h,a^{1},a^{2},a^{3})\right],
\label{darbouxtransformcyclide}
\end{equation}
 where $V(h,a^{1},a^{2},a^{3}):\R^{2}\to \mathcal{L}$ is the
 solution of the $D_{h}$-linear system for $\Phi=1$
 with initial condition $X(a^{1},a^{2},a^{3})$ given by 
(\ref{transformingfunctions}).
Thus,
 $D(h,a^{1},a^{2},a^{3})(f_{m})$ is a special isothermic immersion with
 Calapso potential $\Phi(h,a^{1},a^{2},a^{3})$.

\end{ex}

\section{The superposition formula and two-soliton solutions}\label{se.superpo}

\begin{thm} Let $\Phi:U\to \R$ be a wave potential with character
$s$, and let $V=(v^{0},...,v^{4})^t$ and $W=(w^{0},...,w^{4})^t$ be
two systems of $h$ and $k$-transforming functions, respectively, $h\neq k$. 
Let $\Phi_{1}=E(\Phi,V)$ and $\Phi_{2}=E(\Phi,W)$ be the corresponding
B\"acklund transforms. Then,
\begin{equation} \Phi_3= \Phi - (h-k)\frac{{v}^3{w}^3 }{\langle V,W\rangle}
\left[\frac{\Phi_2 -
\Phi_1}{\Phi_1\Phi_2}\right].\label{sovrapposizione}
\end{equation}
is a wave potential with character $s$.
\end{thm}

\begin{figure}[ht]
\begin{center}
\fbox{\includegraphics[height=5.5cm,width=7cm]{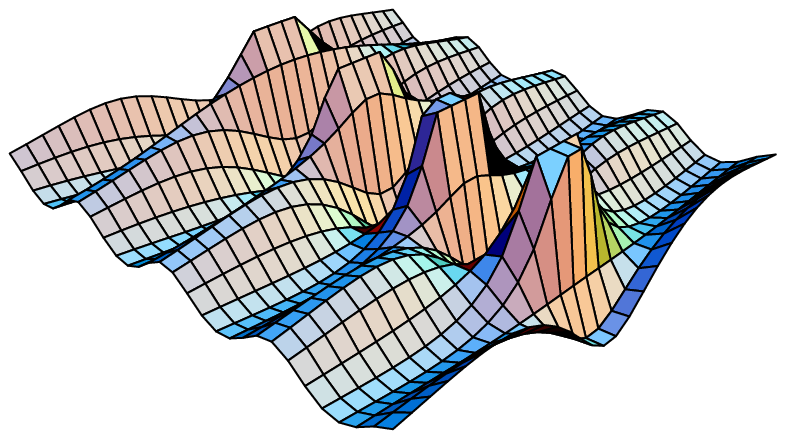}\vrule
\includegraphics[height=5.5cm,width=7cm]{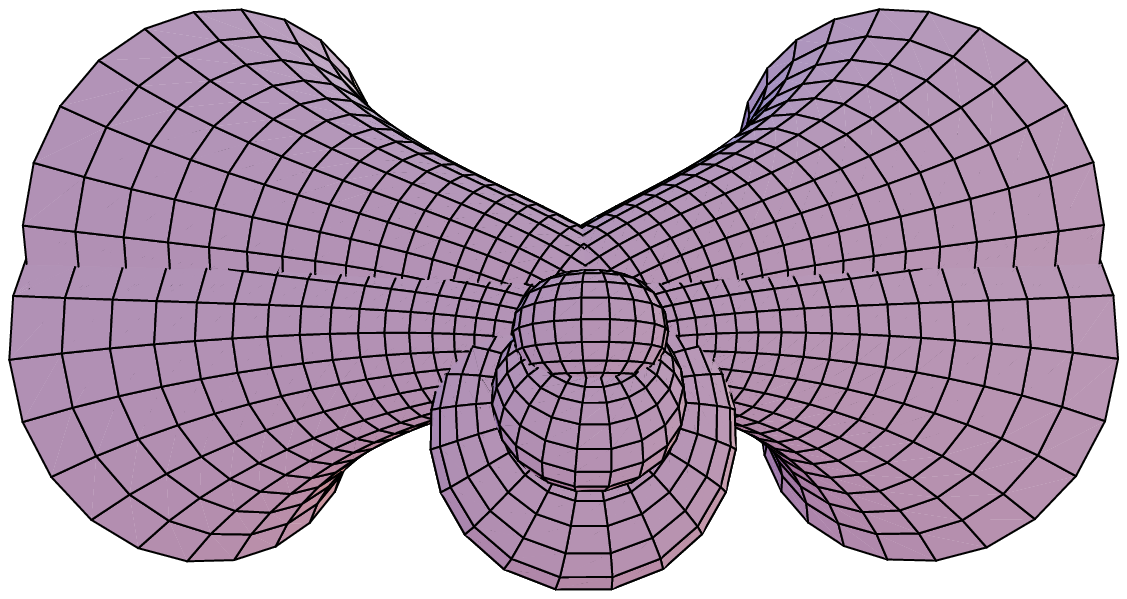}}
\caption{the superposition of the wave potentials
$\Phi(4/3,0,1/100,1/300)$ and $\Phi(10/9,0,1/500,0)$ and the
corresponding special isothermic surface with $m=1$.}
\end{center}
 \end{figure}

\noindent \begin{proof}{\rm Let $B(m):U\to G$ be the normal frame
with potential $\Phi$ and spectral parameter $m$, $m\neq h$.
Consider the curved flat framing $A:=B(m)g^{+}(V)$ and its dual
framing $\bar{A}:=AJ$; then $\bar{A}$ is a curved flat framing and
the corresponding normal frame $\bar{B}(m)$ is computed to be $$
\bar{B}(m) =\bar{A}g^{+}(T)=B(m)g^{+}(V)Jg^{+}(T), $$ where
$T:U\to \mathcal{L}$ is the smooth map defined as in (\ref{Tequation}).
Next, consider the map $Y:U\to \mathcal{L}$ given by
\begin{equation}
\left(\frac{(h-k)w^{4}}{(h-m)v^{4}},\frac{v^{1}w^{4}-v^{4}w^{1}}{v^{4}},
\frac{v^{4}w^{2}-v^{2}w^{4}}{v^{4}},\frac{v^{4}w^{3}-v^{3}w^{4}}{v^{4}},
-\frac{(h-m)}{(h-k)}\langle {V},{W}\rangle\right)^t.
\label{Yequation}
\end{equation}
Then $\tilde{A}:= \bar{A}g^{+}(Y)$ is a curved flat framing with
spectral parameter $k-m$ and Calapso potential ${\Phi_1}$.
Further, let $L:U\to \mathcal{L}$ be defined by
\begin{equation}
{L}={T}^{-1}\star{Y}=\left(-\langle{T},{Y}\rangle,
\frac{y^{1}t^{4}-y^{4}t^{1}}{t^{4}},\frac{y^{2}t^{4}-y^{4}t^{2}}{t^{4}},
\frac{y^{3}t^{4}-y^{4}t^{3}}{t^{4}},\frac{y^{4}}{t^{4}}\right)^t,
\label{Lequation}
\end{equation}
so that $$
\tilde{A}=[\bar{A}g^{+}(T)][g^{+}(T^{-1})g^{+}(Y)]=\bar{B}(m)g^{+}(L).
$$ Since $\bar{B}(m)g^{+}(L)$ is a curved flat framing with
Calapso potential $\Phi_{1}$ and spectral parameter $k-m$, then
$L$ is a solution of the $D_{k}$-system with potential $\Phi_{1}$.
Combining (\ref{Tequation}), (\ref{Yequation}), (\ref{Lequation})
and using the constraints $C_{2}(V)=C_{2}(W)=0$, it is a
computational matter to check that also $L$ satisfies the
constraint $$
C_{2}(L)=\Phi_{1}^{2}L^{0}-kL^{3}+\frac{s}{2}\Phi_{1}^{-2}L^{4}=0.
$$ This implies that $L$ is a system of $D_{k}$-transforming
functions for the potential ${\Phi_1}$. In particular, the
B\"acklund transform $E(\Phi_{1},L)$ is a new solution of the
differential equation (\ref{solitoneq}). On the other hand,
$E(\Phi_{1},L)$ is computed to be
\begin{equation}
E(\Phi_{1},L)=\frac{\left(
v^{1}w^{1}+v^{2}w^{2}+v^{3}w^{3}-v^{0}w^{4}-v^{4}w^{0}\right)\Phi^{2}-(h-k)
\left(w^{3}v^{4}-w^{4}v^{3}\right) }{\left(
v^{1}w^{1}+v^{2}w^{2}+v^{3}w^{3}-v^{0}w^{4}-v^{4}w^{0}\right)\Phi},
\end{equation}
from which follows that $$ E(\Phi_{1},L) = \Phi -
{v}^3{w}^3\frac{(h-k) }{\langle V,W\rangle} \left[\frac{\Phi_2 -
\Phi_1}{\Phi_1\Phi_2}\right]. $$ }
\end{proof}

\begin{remark}
The wave potential $\Phi_{3}$ can be realized as
the Calapso potential of the special isothermic immersion
\begin{equation}
[\left(B(m)g^{+}(V)Jg^{+}(Y)\right)_4],\label{foursurface}
\end{equation}
 which can be explicitly computed from the normal frame
$B(m)$ and the two sets $V$ and $W$ of transforming functions. 
\end{remark}

\begin{ex}
[\textbf{Two-solitons and the corresponding isothermic surfaces}]
Two-soliton solutions can be computed by means of
(\ref{transformingfunctions}) and (\ref{sovrapposizione}). We then
obtain a six-parameter family of wave potentials given by $$ \Phi
= 1 - \frac{(h-k)(w^{3}v^{4}-w^{4}v^{3})}{\langle V, W\rangle}, $$
where $V$ and $W$ are the $h$ and $k$-transforming functions
corresponding to the initial conditions $X(a^{1},a^{2},a^{3})$ and
$X(b^{1},b^{2},b^{3})$.
Special isothermic immersions with spectral parameter $m$ and
Calapso potentials $\Phi$ can be constructed by using
(\ref{foursurface}), (\ref{transformingfunctions}) and the
explicit formulae for the central frame fields $B(m)$ of Dupin
cyclides.
\end{ex}

\bibliographystyle{amsalpha}

\end{document}